\documentclass[11pt]{amsart}
\usepackage{amsmath}

\def\ol{\overline}
\def\lf{\left}
\def\ri{\right}
\def\a{{\alpha}}
\def\b{{\beta}}
\def\g{{\gamma}}
\def\d{{\delta}}
\def\x{{\xi}}
\def\z{{\zeta}}

\def\bb{{\bar\b}}
\def\bd{{\bar\d}}

\def\bz{{\bar\zeta}}
\def\abb{{\a\bb}}
\def\gbd{{\gamma\bd}}
\def\xbz{{\x\bz}}
\def\boldx{{\bold x}}
\def\wt{\widetilde}
\def\tn{{\wt\nabla}}
\def\p{\partial}
\def\re{{\mathbb R}}
\newcommand\ce{{\mathbb C}}

\def\R{R_\abb}

\def\G{g_{ij}}

\def\r{R_{ij}}

\def\K{K\"ahler }
\def\KR{K\"ahler-Ricci }
\def\KRF{K\"ahler-Ricci flow }
\def\KRS{K\"ahler-Ricci soliton }
\def\KRS{K\"ahler-Ricci soliton }

\def\be{\begin{equation}}
\def\ee{\end{equation}}
\newtheorem{thm}{Theorem}[section]
\newtheorem{lem}{Lemma}[section]

\theoremstyle{definition}

\theoremstyle{remark}

\numberwithin{equation}{section}

\begin{document}

\title{gradient K\"ahler-Ricci solitons and a uniformization conjecture}

\author{Albert Chau$^{1}$}
\address{Harvard University, Department of Mathematics,
  One Oxford Street, Cambridge, MA 02138, USA}
\email{chau@math.harvard.edu}

\author{Luen-Fai Tam$^2$}
\thanks{$^{1}$Research
partially supported by The Institute of Mathematical Sciences, The Chinese
University of Hong Kong, Shatin, Hong Kong, China.}
\thanks{$^2$Research
partially supported by Earmarked Grant of Hong Kong \#CUHK4032/02P}

\address{Department of Mathematics, The Chinese University of Hong Kong,
Shatin, Hong Kong, China.} \email{lftam@math.cuhk.edu.hk}

\renewcommand{\subjclassname}{%
  \textup{2000} Mathematics Subject Classification}

\date{October 2003.}



\begin{abstract}
 In this article we study the limiting behavior of the \KRF on complete non-compact \K
 manifolds.  We provide sufficient conditions under
 which a complete non-compact gradient K\"ahler-Ricci soliton is
biholomorphic to $\ce^n$.  We also discuss the uniformization conjecture by Yau
\cite{Y} for complete non-compact \K manifolds with positive holomorphic bisectional
curvature.
\end{abstract}

\maketitle

\markboth{Albert Chau and Luen-Fai Tam} {gradient K\"ahler-Ricci solitons and a
uniformization conjecture}
\section{introduction}\label{intro}
In this paper, we show when a complete non-compact gradient \KRS is biholomorphic
to $\ce^{n}$.  We will also discuss when a general solution to the \KRF on a
non-compact \K manifold converges after rescaling to a complete flat \K limit
metric.

Canonical examples of such solitons on $\ce^{n}$ were first provided by
Cao\cite{cao, Cao}.  These examples are all rotationally symmetric with positive
holomorphic bisectional curvature.  It would be interesting to know how many
other complete gradient \KRS metrics there are on $\ce^{n}$.  Our results may be
of use here. Another reason for our interest in gradient \K Ricci solitons is
that they may serve as models for the uniformization conjecture by Greene-Wu \cite{GW1},  Siu \cite{Siu} and in the most general form  by Yau \cite{Y}
which states that any complete non-compact \K manifold with positive holomorphic
bisectional curvature is biholomorphic to $\ce^{n}$.  Using our techniques and
ideas we shed light on recent approaches to proving this conjecture using the
\KRF \cite{Sh,Sh2, NT}.

 A gradient Ricci soliton is defined as follows. Let $\G(x, t)$ be a family of metrics on a Riemannian manifold
$M$ satisfying the Ricci flow equation:
\begin{equation}\label{rf}
\frac{\p}{\p t}\G=-2\r-2\rho \G.
\end{equation}
for $0\le t<\infty$, where $\r$ denotes the Ricci tensor at time $t$ and $\rho$
is a constant.   $\G(x, t)$ is said to be a {\it gradient Ricci soliton of steady
type}, if $\rho=0$ and if there is a potential function $f$ and a family of
diffeomorphisms $\varphi_t$ generated by the gradient of   $-f$  with respect to
$\G(x, 0)$ such that $\G(x, t)=\varphi_t^*\lf(\G(x, 0)\ri)$. If $\rho>0$
(respectively $\rho<0$), then it is said to be of {\it expanding type}
(respectively {\it shrinking type}). If $\G(x, t)$ is a gradient Ricci soliton
with potential function $f$ then one has
\begin{equation}\label{rs}
f_{ij}=2R_{ij}(x,0)+2\rho g_{ij}(x,0)
\end{equation}
where $f_{ij}$ is the Hessian of $f$ with respect to $\G(x, 0)$.

If $(M, g_\abb(x, 0))$ is a \K manifold,   (\ref{rf}) is referred to as the \KRF and
is written as
\begin{equation}\label{krf}
\frac{\p}{\p t} g_\abb=-\R-2\rho  g_\abb.
\end{equation}
A gradient Ricci-soliton solution to (\ref{krf}) is referred to as a gradient
\KRS.  In this case (\ref{rs}) takes the form
\begin{equation}\label{krs}
\begin{split}
f_{\abb}&=\R+2\rho  g_\abb\\
f_{\a\b}&=0.
\end{split}
\end{equation}
Hence the gradient of $f$ is a holomorphic vector field and the diffeomorphism
$\varphi_t$ is a biholomorphism.  At times,  we may refer to a Riemannian
manifold $(M, \G)$ as a Ricci-soliton if the corresponding solution to (\ref{rf})
is a Ricci-soliton .  We do likewise in the \K case.

   We consider gradient \KR solitons which are either (i) steady with positive Ricci curvature so that the scalar curvature attains maximum at some point; or (ii) expanding with nonnegative Ricci curvature. Under either of these conditions, it is not hard to prove that there is a unique equilibrium point $p$ where the gradient of the potential function $f$ is zero. Our main result for gradient \K Ricci solitons is:

\begin{thm}\label{mainthm} Let $(M, g_\abb)$ be a complete non-compact gradient \KRS with potential $f$ satisfying  either of the conditions mentioned above,   and
let $g_\abb(x, t)$ be the corresponding solution to (\ref{krf}).  Let $p$ be the
equilibrium point and let $\bold v_p \in T^{1,0}_p (M)$ be a fixed nonzero vector
with $|\bold v_p|_0=1$.   Then for any sequence of times ${t_k}\rightarrow
\infty$, the sequence of complete \K metrics $ \frac{1}{ |\bold v_p|_{t_k}^{2} }
g_\abb(x, t_k)$ subconverges on compact sets of $M$ to a complete flat \K metric
$h_{\abb}$ on $M$ if and only if  $ R_\abb(p)=\beta g_\abb(p)$ at $t=0$ for some
constant   $\beta$.  In particular, if the condition is satisfied then $M$ is
biholomorphic to $\ce ^n$.
\end{thm}
Here for a tangent vector $\bold v$ on $M$, $|\bold v|_t$ denotes the length of
$\bold v$ in the metric $g(t)$.

Next, we consider   general complete non-compact K\"ahler manifolds with nonnegative
holomorphic bisectional curvature.  In \cite{Sh, Sh2} (see also \cite{NT}), W.-X.
Shi proved that on a complete noncompact K\"ahler manifold $(M,g_\abb)$ with bounded
nonnegative holomorphic  bisectional curvature such that
\begin{equation}\label{scalar}
\frac{1}{V_x(r)}\int_{B_x(r)}R \le \frac{C}{1+r^2}
\end{equation}
for some constant $C$ for all $x\in M$ and for all $r$, the K\"ahler-Ricci flow
$$
\frac{\p}{\p t}g_\abb(x,t)=-R_\abb(x,t)
$$
with initial condition $g_\abb(x,0)=g_\abb(x)$ has a long time solution. Moreover,
useful estimates were obtained. In \cite{Sh}, an approach by Shi to prove the
uniformization conjecture of Greene-Wu-Siu-Yau for manifolds satisfying
(\ref{scalar}) is to use the K\"ahler-Ricci flow  to produce a complete flat \K
metric $h$ the   \K manifold $M$. More precisely, one considers the rescaled metrics
$\frac1{|\mathbf v_p|^2_t}g_\abb(x,t)$ and shows that a subsequence will converge to
a flat complete \K metric $h$. Here $\mathbf v_p$ is a fixed vector in
$T_p^{1,0}(M)$ and $|\mathbf v_p|_t$ is its length in $g(t)$. However, the proof in
\cite{Sh} is not quite satisfactory. First, as noted in \cite{Chen} the completeness
of $h$ is unclear from \cite{Sh} and has yet to be verified.  On the existence of
$h$, the authors would like to point out that the proof in \cite{Sh} depends
critically on a bound for a  quantity $Q$ (see (\ref{qbound}) for more details) and
that Shi's proof of this bound appears to be incorrect. More specifically, the
formula on \cite[p.156]{Sh} for $\frac{\p}{\p t}Q$ seems to be incorrect. In this
paper we partially rectify these issues by providing a proof for the completeness of
$h$ assuming we have an a priori bound for $Q$.  We do this in section 4
(Theorem{\ref{firstway}). In general, in the absence of such a bound, we prove that
completeness is in many cases a natural condition that follows from the existence of
$h$ alone.  In this direction our main result is:
\begin{thm}\label{secondway1} There exists a constant $C(n)$ depending only on $n$ such that if $M^n$ is a complete noncompact K\"ahler manifold with bounded nonnegative holomorphic bisectional curvature
satisfying:
\begin{enumerate}
\item[(i)]
 $$
\frac1{V_x(r)}\int_{B_x(r)}R\le \frac{C(n)}{1+r^2}
$$
for all $x\in M$ and for all $r>0$; and
\item[(ii)]   there exist  a point $p \in M$ and a sequence $t_k\to\infty$ such that at $p$   $\frac{1}{|\mathbf v_p|^2_{t_k}}g(p,t_k)$ are uniformly equivalent to $g(p,0)$, where $\mathbf v_p$ is a fixed vector in $T_p^{1,0}(M)$ with $|\mathbf v_p|_0=1$.
\end{enumerate}  Then the metrics $\frac{1}{|\mathbf v_p|^2_{t_k}}g(x,t_k)$
subconverge uniformly in the $C^\infty$ topology in compact sets to a complete K\"ahler flat metric
on $M$. In particular, the universal covering space of $M$ is biholomorphic to
$\mathbb C^n$.
\end{thm}

The authors would like to thank Prof. S.T.Yau for helpful discussions and
support. The first Author would also like to thank Prof. Richard Hamilton for
helpful discussions.

\section{A necessary condition for convergence}
In this section, we prove the necessary part of Theorem \ref{mainthm}. In fact,
we have the following:

\begin{thm} Let $\G(x,t)$ be a gradient Ricci soliton with function $f$
and diffeomorphisms $\varphi_t$ generated by  $\nabla_0 (-f)$, where $\nabla_0$
is the covariant derivative respect to $\G(x, 0)$. Suppose the flow $\varphi_t$
has an equilibrium point $p$ and suppose there exist a subsequence $t_k\to\infty$
and positive numbers $\sigma(t_k)$ such that $\sigma(t_k)\G(x, t_k)$ converges
uniformly in a neighborhood of $p$ to a Riemannian metric $h_{ij}$. Then at
$t=0$, $\r(p) =\beta \G(p) $  for some constant $\beta$.
\end{thm}

\begin{proof} In the following $\G(x,0)$ will simply be denoted by $g$ and the metric at
time $t$ will be denoted explicitly by $g(t)$.

Choose a coordinate neighborhood $V$ of $p$ with coordinates $\mathbf
x=(x^1,\dots, x^n)$ such that $\mathbf x(p)=0$, $g_{ij}(0)=\delta_{ij}$,
$\frac{\p}{\p x^k}g_{ij}(0)=0$, $\frac{\p^2f}{\p x^i\p
x^j}(0)=\lambda_i\delta_{ij}$. By (1.2),   it is sufficient to prove that
$\lambda_{i}=\lambda_j$ for all $i,\ j$.

Let $\mathbf  v_0\in T_p(M)$ such that
\begin{equation}
\mathbf  v_0=\sum_kv_0^k\frac{\p}{\p x^k}.
\end{equation}
Let $F^i(\mathbf  x)=g^{ij}(\mathbf x)\frac{\p f}{\p x^j}(\mathbf x).$  Since $p$
is an equilibrium point, $\varphi_t(0)=0$ for all $t$, $\nabla_0 f(p)=0$. Hence
$F^i(0)=0$ and $\frac{\p F^i}{\p x^j}(0)=\lambda_i\delta_{ij}$.

We may assume that there is a constant $C_1$ such that $|F(\mathbf  x)|\le
C_1|\mathbf  x|$ on $V$, where $|\mathbf  x|^2=\sum_i(x^i)^2$. Hence for any
$T>0$, there exists a constant $a>0$ such that the equation
\begin{equation}\begin{cases} \frac{d\bold x}{dt}=&F(\bold x(t))\\
\boldx(0)=&\bold x_0
\end{cases}
  \end{equation}
has a unique solution on $[0,T]$ with image inside $V$ whenever $|\mathbf
x_0|^2=\sum_i\lf(x_0^i\ri)^2 \le a^2$.

Consider the curve $\a(s)= (sv^1_0,\dots,sv^n_0)$ so that $\a'(0)=\mathbf  v_0$.
There exists $s_0>0$ such that $|\a(s)|\le a$ for all $0\le s\le s_0$. Hence for
all $0\le s\le s_0$, the solution $\mathbf x\lf(t;\a(s)\ri)$ of (2.2) with
initial value $\a(s)$ is defined on $0\le t\le T$ with image inside $V$. Since
$\varphi_t(0)=0$ for all $t$, $\lf(\varphi_t\ri)_*(\mathbf  v_0)=\frac{\p}{\p
s}\varphi_t(\a(s))\Big|_{s=0}\in T_p(M)$. Denote $\lf(\varphi_t\ri)_*(\mathbf
v_0)$ by
\begin{equation}
\sum_kv^k(t)\frac{\p}{\p x^k}.
\end{equation}
In local coordinates $\varphi_t(\alpha(s))=\mathbf
x(t;\a(s))=\lf(x^1\lf(t;\a(s)\ri),\dots,x^n\lf(t;\a(s)\ri)\ri)$ and
\begin{equation}
\frac{\p}{\p s}\varphi_t(\a(s))\Big|_{s=0}=\frac{\p x^k}{\p
s}\lf(t;\a(s)\ri)\Big|_{s=0}\frac{\p}{\p x^k}.
\end{equation}
Hence $v^k(t)$ is given by
\begin{equation}
v^k(t)=\frac{\p x^k}{\p s}\lf(t;\a(s)\ri)\Big|_{s=0}.
\end{equation}
Now for $0\le t\le T$,
\begin{equation}
\begin{split}
 \frac{d}{dt}v^k(t) &=\frac{\p}{\p t}\frac{\p x^k}{\p s}\lf(t;\a(s)\ri)\big|_{s=0}\\
&=\frac{\p}{\p s} \frac{\p x^k}{\p t}\lf(t;\a(s)\ri)\big|_{s=0}\\
&=\frac{\p}{\p s}F^k\lf(\mathbf x(t;\a(s))\ri)\big|_{s=0}\\
&=\lf[\frac{\p}{\p x^i}F^k\lf(\mathbf x(t;\a(s))\ri)\frac{\p x^i}{\p s}\lf(\mathbf x\lf(t;\a(s)\ri)\ri)\ri]_{s=0} \\
&=\lambda_{k}\delta_{ki}v^i(t) \\
&=\lambda_{k}v^k(t),
\end{split}
\end{equation}
where we have used (2.5), the fact that $\mathbf x(t;\a(0))=\mathbf x(t;0)=0$
because $F(0)=0$, and the fact that $\frac{\p F^i}{\p x^k}=\lambda_i\d_{ik}$ at
$0$. Using the initial condition, we conclude that
\begin{equation}
v^k(T)=\exp(\lambda_kT)v^k_0
\end{equation}
and
\begin{equation}
\lf(\varphi_T\ri)_*\lf(\sum_k   v_0^k\frac{\p}{\p x^k}\ri)=
\sum_k\exp(\lambda_kT)v^k_0\frac{\p}{\p x^k}.
\end{equation}
Hence
\begin{equation}
|\mathbf  v_0|^2_{g(T)}=|\mathbf  v_0|^2_{\varphi_T^*(g_0)}
=|\lf(\varphi_T\ri)_*\mathbf  v_0|^2_{g}=\sum_i\exp(2\lambda_iT)\lf(v^i_0\ri)^2
\end{equation}
for all $T>0$. Suppose there exist $t_k\to\infty$, $\sigma(t_k)g(t_k)$ converges
in $C^\infty$ topology to a Riemannian metric $h $ on a neighborhood of $p$. Then
there exists a constant $C_2>0$ such that for any $\mathbf  v,\ \mathbf  w\in
T_p(M)$ with $|\mathbf  v|_{g}=|\mathbf  w|_{g}$, and for all k, we have
\begin{equation}
C_2^{-1}\le \frac{|\mathbf  v|^2_{\sigma(t_k)g(t_k)}}{|\mathbf
w|^2_{\sigma(t_k)g(t_k)}}\le C_2.
\end{equation}
In the coordinates $(x^1,\cdots,x^n)$, by (2.9), we have
\begin{equation}
C_2^{-1}\le \frac{\sum_i\exp(2\lambda_it_k)\lf(v^i_0\ri)^2
}{\sum_i\exp(2\lambda_it_k)\lf(w^i_0\ri)^2 }\le C_2
\end{equation}
for all $k$, whenever $\sum_i \lf(v^i_0\ri)^2=\sum_i \lf(w^i_0\ri)^2=1$. Since
$t_k\to\infty$,  $\lambda_i=\lambda_j$ for all $i$ and $j$.
\end{proof}

\section{a sufficient condition for convergence}\label{maintheorem}
In this section we prove the sufficient part of Theorem \ref{mainthm}. First, we
have the following on the existence of an equilibrium point.

\begin{lem}\label{fixedpoint} Let $(M^n, g_\abb)$ be a complete non-compact gradient \KRS with potential $f$
satisfying either of the following two conditions:
\begin{enumerate}
\item At $t=0$, $f_\abb=R_\abb$   and $R_\abb > 0$ so that the scalar curvature
$R$ attains maximum at some point in $M$. \item At $t=0$, $f_\abb= R_\abb+g_\abb$
and $ R_\abb \ge 0$.
\end{enumerate}  Then there is a unique point $p \in
M$ at which $\nabla_0 f(p)=0$, where $\nabla_0$ is the covariant derivative with
respect to $g(0)$. Also, $M$ is diffeomorphic to $\re^{2n}$.
\end{lem}

\begin{proof} It will suffice to show that $f$ is a strictly convex exhaustion function, see \cite[Theorem 3]{GW}.

In case (1), this follows from  the proof of  \cite[Theorem 20.1]{Ha4}, see
also\cite{ CH}.

 In case (2), we begin by noting that (2) together
with (\ref{krs}) imply that the Hessian of $f$ with respect to $g(0)$ satisfies $
D^2f \ge g (0)$, thus $f$ is indeed strictly convex. Next, let $q$ be a fixed
point and  consider an arbitrary geodesic $\gamma(s)$ originating at $q$
parametrized by arc length in $g(0)$. Then along $\gamma(s)$ we have
$\frac{d^2f}{ds^2}(\gamma(s))=D^2 f(\gamma '(s), \gamma '(s))\geq 1$. Integrating
this we get
\begin{equation}
\begin{split}
f(\gamma(s))-f(q) =&
f(\gamma(s))-f(\gamma(0))\\
  =&\int^{s}_{0} \lf[\int^{\tau}_{0}\frac{d^2
f}{ds^2}(\gamma(\rho))d\rho - \frac{df}{ds}(\gamma(0))\ri]d \tau \\
\geq &\int^{s}_{0} \lf(\int^{\tau}_{0}d\rho - |\nabla_0 f|(q)\ri)d \tau\\
\geq &\int^{s}_{0} (\tau - |\nabla_0 f|(q))d \tau\\
\geq &\frac{s^2}{2}-|\nabla_0 f|(q)s
\end{split}
\end{equation}
where $\nabla_0$ is the covariant derivative with respect to $g(0)$.
  It is now
clear that $f$ is an exhaustion function on $M$.  This completes the proof of the
lemma.
\end{proof}

The sufficient part of Theorem  \ref{mainthm}  will follow from Lemma
\ref{fixedpoint} and the following lemmas.  In the following, when we say   case
(1) (respectively case (2)), we mean that the potential $f$ in Theorem
\ref{mainthm}  satisfies condition (1) (respectively condition (2)) in Lemma
\ref{fixedpoint}.

Let $p$ be the equilibrium point in Theorem  \ref{mainthm}, whose existence is
implied by Lemma \ref{fixedpoint}. $B_t(R)$ will denote the geodesic ball of
radius $R$ with respect to the metric $ g (t)$ with center $p$. In particular
$B_0(R)$ is the geodesic ball of radius $R$ with respect to the initial metric
$g(0)$.

\begin{lem}\label{shrink} With the same assumptions and notations as in Lemma \ref{fixedpoint}, for any $R>0$, the following are true:
\begin{enumerate}
\item[{(i)}] $B_{t_1}(R)\subset B_{t_2}(R)$ for all $t_1\le t_2$; \item[{(ii)}]
for any $T\ge0$, $q\in B_T(R)$, $\bold w_q\in T^{1,0}(M)$,
$$
|\bold w_q|_t\le  \exp\lf(-C_R(t-T)\ri)|\bold w_q|_T
$$
for all $t\ge T$, where    $C_R>0$ is a constant depending only on $R$ and
$g(0)$; and \item[{(iii)}] for any integer $k\ge 0$, for any $t\ge0$,
$$
||\nabla_t^k Rm(t) ||_{g(t)}\le C(R,k)
$$
on $B_0(R)$ for some constant $C(R,k)$ depending only on $R$, $k$ and $g(0)$,
where $\nabla_t$ is the covariant derivative with respect to $g(t)$ and $Rm(t)$
is the curvature tensor of $g(t)$.
\end{enumerate}
\end{lem}
\begin{proof} Let $\varphi_t$ be the biholomorphism of $M$ generated by the gradient of $-f$ so that $g(t)=\varphi^*_t(g(0))$. Then  $\varphi_t(p)=p$ by the definition of $p$. Since $R_\abb\ge0$ in both cases in the assumptions of Lemma \ref{fixedpoint}, $g_\abb(t_2)\le g_\abb(t_1)$ if $t_1\le t_2$. From these, it is easy to see that (i) is true.

Since $\varphi_t:(M,g(t))\to  (M,g(0))  $ is an isometry and $\varphi_t(p)=p$,
$\varphi_t$ will map  $(B_t(R),g(t))$  isometrically onto $(B_0(R),g(0))$. Hence
by (i) if $t\ge T$, the greatest lower  bound of the Ricci curvature of $g(t)$ in
$B_T(R)$ is no less than the greatest lower bound of the Ricci curvature of
$g(T)$ in $B_T(R)$, which is the same as the greatest lower bound of the Ricci
curvature of $g(0)$ on $B_0(R)$.

Now let $q\in B_T(R)$ and if $\bold w=\bold w_q\in T^{1,0}(M)$,
\begin{equation}
\begin{split}
\frac{\p}{\p t}\lf(g_\abb(q,t)w^\a w^{\bb}\ri)&=\lf(-R_\abb(q,t)-2\rho g_\abb(q,t)\ri)w^\a w^\bb\\
&\le -C_1 g_{\abb}(q,t) w^a w^\bb
\end{split}
\end{equation}
for some constant $C_1>0$ depending only on $R$ and $g(0)$. In fact, if case (2)
is assumed so that $\rho=1/2$, then $C_1$ can be taken to be 1. If case (1) is
assumed so that $\rho=0$ then $C_1$ can be taken to be twice the greatest lower
bound of the Ricci curvature of $g(0)$ in $B_0(R)$,  which is positive. Dividing
both sides of the above inequality by $g_{\abb}(q,t) w^a w^\bb$ and  integrating
from $T$ to $t$, (ii) follows.

Since $B_0(R)\subset B_t(R)$ for $t\ge0$ and since $(B_t(R),g(t))$ is isometric
to $(B_0(R),g(0))$, it is easy to see that (iii) is true.
\end{proof}

\begin{lem}\label{uniformshrink} With the same assumptions and notations as in Lemma \ref{fixedpoint}, let $R>0$ and $T\ge0$.
Then there exists a constant $C_R>0$ which depends only on $R$ and $g(0)$ with
the following property:  For any $q\in B_{T}(R)$, $\bold u_p\in T^{1,0}(M)$,
$\bold w_q\in T^{1,0}(M)$ with $|\bold u_p|_{T}=|\bold w_q|_{T}$,
\begin{equation}
C_R^{-1}\le \frac{|\bold u_p|_{t}}{|\bold w_q|_t}\le C_R \end{equation} for all
$t\ge T$.
\end{lem}
\begin{proof} For any $t\ge T\ge0$, let $q$, $\bold w_q$ and $\bold u_p$ as in the assumptions. Let $\gamma_t(s)$ be a minimal geodesic from $p$ to $q$ in the metric $g(t)$. Let $\bold w(s)$ be a parallel vector field with respect to $g(t)$ along $\gamma_t(s)$ such that $\bold w(\gamma_t
(d))=\bold w_q$, where $d=d_t(p,q)$ is the distance between $p$ and $q$ in
$g(t)$. Then in both case (1) and case (2)
\begin{equation}\label{estimatericiatq}
\begin{split}
\frac{\p}{\p t}\log \lf[\frac{g_\abb(p,t)u_p^\a u_p^\bb}{g_\abb(q,t)w_q^\a w_q^\bb}\ri]&=\frac{R_\abb(q, t)w^\a(d) w^\bb(d)}{ g_\abb(q, t)w^\a(d) w^\bb(d)}-\frac{ R_\abb(p, t)u_p^\a u_p^{\bb}}{g_\abb(p, t)u_p^\a v_p^{\bb}}\\
&= \frac{R_\abb(q, t)w^\a(d) w^\bb(d)}{ g_\abb(q, t)w^\a(d) w^\bb(d)}-\frac{ R_\abb(p, t)w^\a(0) w^{\bb}(0)}{g_\abb(p, t)w^\a(0) w^{\bb}(0)}\\
&= \int_{0}^{d}\frac{\p}{\p s}\lf(\frac{ R_\abb(\gamma_t(s), t)w^\a(s)
w^{\bb}(s)}{ g_\abb(\gamma_t(s), t)w^\a(s)
w^{\bb}(s)}\ri)ds \\
&\leq d\max_{0\le s\le d}||\nabla_t R_\abb||_{g(t)}(\gamma_t(s),t).\end{split}
\end{equation}
Here we have used the assumption that $R_\abb(p,0)=\beta g_\abb(p,0)$ for some
constant $\beta$ and hence $R_\abb(p,t)=\beta g_\abb(p,t)$ for all $t$ because
$\varphi_t(p)=p$. Since $q\in B_T(R)$, by Lemma \ref{shrink}(ii) for all $t\ge
T$, we have
$$
d\le R\exp(-C_1(t-T))
$$
for some positive constant $C_1$ depending only on $R$ and $g(0)$. By Lemma
\ref{shrink}(iii),   we have
\begin{equation}\label{estimatemetricatpq}
 \lf|\frac{\p}{\p t} \log \lf[\frac{g_\abb(p,t)u_p^\a u_p^\bb}{g_\abb(q,t)w_q^\a w_q^\bb}\ri]\ri| \le C_2 R\exp\lf(-C_1(t-T)\ri)
\end{equation}
where   $C_2$ is a constant  depending only on $R$ and $g(0)$. Integrating
(\ref{estimatemetricatpq}) from $T$ to $t$, using that fact that $|\bold
u_p|_T=|\bold w_q|_T$, the resut follows.
\end{proof}
\begin{lem}\label{existenceoflimit}
With the same assumptions as in Lemma \ref{fixedpoint}, for any sequence of times
${t_k}\rightarrow \infty$, the sequence of complete \K metrics $h( k)
=\frac{1}{|\bold v_p|^2_{t_k}} g( t_k)$ has a subsequence converging in
$C^{\infty}$ on compact sets of $M$ to a flat \K metric $H$ on $M$.
\end{lem}
\begin{proof} For any $t\ge 0$, let $h(t)= \frac{1}{\sigma(t)} g( t)$, where $\sigma(t)=|\bold v_p|^2_{t}$. In the following $\widehat{Rm}(t)$ and $\widehat\nabla$ will denote the curvature tensor and the covariant derivative of $h(t)$, and $ {Rm(t)}$ and $ \nabla$ will denote the curvature tensor and the covariant derivative of $g(t)$.

 By Lemma \ref{shrink}, for any interger $m\ge0$ and $R>0$, there is a constant $C_1$ depending only on $m$,  $R$ and $g(0)$ such that  \begin{equation}\label{curvaturedecay}
||\widehat\nabla^m\widehat{Rm}(t)||^2_{h(t)}=
\sigma^{m+2}(t)||\nabla^mRm||_{g(t)}\le C_1\sigma^{m+2}(t)
\end{equation}
on $B_0(R)$ with
\begin{equation}\label{sigma}
\sigma(t)\le
\exp(-C_2t)
\end{equation}
 for some   constant $C_2>0$ depending only on $g(0)$. By Lemma \ref{uniformshrink} and the definition of $h(t)$, there is
a constant $C_3>0$ depending only on $R$ and $g(0)$ such that
\begin{equation}\label{equivalent}
C_3^{-1}g(0)\le h(t)\le C_3g(0)
\end{equation}
for all $t\ge0$.

Let $(z^1,\dots,z^m)$ be a fixed local coordinates in an coordinates neighborhood
$U\subset B_0(R)$. We want to prove that
\begin{equation}\label{dh}
\lf|\frac{\p}{\p z^\x}h_\abb\ri|(x,t)\le C_4
\end{equation}
for some constant $C_4$ for all $x\in U$ and for all $t$.

Let $\Gamma_{\a\b}^\gamma$ be the Christoffel symbols of $h(t)$ which is also the
Christoffel symbols of $g(t)$ in the coordinates $z^\a$ and let
$\wt\Gamma_{\a\x}^\tau$ be the Christoffel symbols of $g(0)$. Let
$A_{\a\x}^\tau=\Gamma_{\a\x}^\tau-\wt\Gamma_{\a\x}^\tau$, then $A_{\a\x}^\tau$ is
a tensor and
$$
 \tn_\x g_\abb= A_{\a\x}^\tau g_{\tau\bb}
$$
where $\tn$ is the covariant derivative with respect to $g(0)$.
 Then the norm of $A$ with respect to $g(0)$ is given by
$$
||A||^2_0=\wt g_{\gamma\bd}\wt g^{\abb} \wt
g^{\x\bz}A_{\a\x}^\gamma\ol{A_{\b\z}^\d}.
$$
By (1.3), we have
\begin{equation}\label{christoffel}
\frac{\p}{\p t}||A||^2_0=-\wt g_{\gamma\bd}\wt g^{\abb} \wt
g^{\x\bz}\lf[g^{\gamma\bar\sigma}\nabla_\a R_{\x\bar\sigma}\ol{A_{\b\z}^\d}
+A_{\a\x}^\gamma \ol{g^{\d\bar\sigma}\nabla_\beta R_{\z\bar\sigma}}\ri]
\end{equation}
Since the equality does not depends on   coordinates, we choose  holomorphic
coordinates $(u^1,\dots,u^n)$ in $U$ such that $\wt g_\abb=\delta_{\a\b}$,
$g_\abb=\lambda_\a\delta_{\a\b}$  at a point. Then
\begin{equation}
\begin{split}
\lf|\wt g_{\gamma\bd}\wt g^{\abb} \wt g^{\x\bz}  g^{\gamma\bar\sigma}\nabla_\a
R_{\x\bar\sigma}\ol{A_{\b\z}^\d}\ri|
&\le C_5 \lf(\sum_{\a,\xi,\lambda}\lambda_\gamma^{-1}\lf|\nabla_\a R_{\x\bar\gamma}\ri|\ri)\,||A||_0\\
&\le C_6 \lf(\sum_{\a,\xi,\lambda}\lambda_\gamma^{\frac 12}\lambda_\a^{-\frac12}\lambda_\xi^{-\frac 12}\lambda_\gamma^{-\frac12}\lf|\nabla_\a R_{\x\bar\gamma}\ri|\ri)\,||A||_0\\
&\le C_7 \exp(-C_7t)||\nabla_\a R_{\x\bar\gamma}||\,||A||_0\\
&\le C_8\exp(-C_7t)||A||_0
\end{split}
\end{equation}
for some constants $C_5-C_8$ depending only on $R$ and $g(0)$ where we have used
Lemma \ref{shrink} and Lemma \ref{uniformshrink}.
 Combining this with (\ref{christoffel}), we have
$$
 \frac{\p}{\p t}||A||_0^2 \le  C_9\exp(-C_7t)||A||_0$$
for some constant $C_9$ depending only on $R$ and $g(0)$.   Since $A=0$ at $t=0$,
we conclude that
$$
||A||_0^2(x,t)\le C_{10}
$$
for some constant $C_{10}$ for all $x\in U$ and for all $t$. From this and (\ref{equivalent}), it is
easy to see that (\ref{dh}) is true.

Now
\begin{equation}\label{ddh1}
\widehat R_{\abb\gamma\bd}=-\frac{\p^2}{\p z^\gamma \p\bar z^{\delta}}h_\abb+h^{\sigma\bar\tau}\lf(\frac{\p}{\p \bar z^\delta}h_{\sigma\bar\beta}\ri)\lf(\frac{\p}{\p z^\gamma}h_{\a\bar\tau}\ri).
\end{equation}
By (3.6)--(3.8), there is a constant $C_{11}$ independent of $t$ such that
$$
 \lf|\Delta_0 h_\abb\ri|=4\lf|\sum_{\gamma}\frac{\p^2}{\p z^\gamma \p\bar z^{\gamma}}h_\abb\ri|\le C_{11}
$$
in $U$. By \cite[Theorem 8.32]{GT}, for any open set $U'\subset\subset U$ there are  constants $C_{12}>0$ and $1>\alpha>0$ independent of  $t$ such that the $C^{1,\a}$ normed of $h_{\abb}$ satisfies
\begin{equation}\label{ddh2}
|h_{\abb}|_{1,\a,U'}\le C_{12}.
\end{equation}
Also, by (3.6)--(3.8), we conclude that
$$
\lf|\frac{\p}{\p z^\sigma}\widehat R_{\abb\gamma\bd}\ri|\le C_{13}
$$
in $U$ for some constant $C_{13}$ independent of $t$.
Hence we can conclude from (\ref{ddh1}) and (\ref{ddh2}) that the $C^\a$ norm of $\Delta_0h_\abb$ in $U'$ is also bounded by a constant independent of $t$. Therefore the $C^{2,\a}$ norm  of $h_\abb$ in any $U'\subset\subset U$ can be bounded by the constant independent of $t$. Similarly, one can prove that the $C^{k,\a}$ norm of $h_\abb$ is bounded by a constant independent of $t$. From this, (\ref{curvaturedecay}), (\ref{equivalent}) and (\ref{sigma}) it is easy to see the lemma is true.
\end{proof}

\begin{lem}\label{completeness} $H$ is complete.
\end{lem}
\begin{proof} We may assume $h(k)$ converge to $H$. Suppose $H$ is not complete. Then there is a
divergent path $\gamma(\tau):[0,\infty)\to M$ from $p$ such that
$\ell_H(\gamma)=L<\infty$, here $\ell_H$ is the length with respect to the metric
$H$. Given $0<\epsilon<L$. Let $a>0$ be such that
$\ell_H(\gamma|_{[0,a]})=L-\epsilon/2$. Since $\gamma|_{[0,a]}$ is compact, there
exists $k_0$ such that for all $k\ge k_0$,
\begin{equation}\label{incompleteL}
L+\epsilon\ge \ell_k(\gamma|_{[0,a]})\ge L-\epsilon
\end{equation}
where $\ell_k$ is the length with respect to $h(t_k)$. By Lemma
\ref{uniformshrink} and by the fact that $h(k)\ge g(t_k)$ because $|\bold
v_p|_{t_k}\le 1$ by Lemma \ref{shrink},  there is a constant which is independent
of $k$ and $k_0$, such that for any $k\ge k_0$, $q\in \wt B_{k_0}(3L)$, $\bold
w_q\in T^{1,0}(M)$, $\bold w_p\in T^{1,0}(M)$ such that if $|\bold
w_q|_{h(k_0)}=|\bold w_p|_{h(k_0)}$, then
\begin{equation}\label{incompleteE}
C^{-1}\le \frac{|\bold w_p|_{h(k)}}{|\bold w_q|_{h(k)}}\le C
\end{equation}
for some   constant $C>0$ depending only on $L$ and $g(0)$. Here $\wt
B_{k_0}(3L)$ is the geodesic ball of radius $3L$ in the metric $h({k_0})$ with
center at $p$.
 Now
reparametrized $\gamma$ by arc length $s$ with respect to  $h(t_{k_0})$. Let
$\gamma|_{(0\le \tau\le a)}=\gamma|_{(0\le s\le b)}$ where $b$ is the length of
$\gamma|_{(0\le \tau\le a)}$ with respect to $h( {k_0})$. By (\ref{incompleteL}),
we have $L-\epsilon\le b\le L+\epsilon$. In particular, $2b\le 3L$ and so
$\gamma([0,2b])\subset\wt B_{k_0}(3L)$. By (\ref{incompleteE}),  for $k\ge k_0$
and for any $b\le s\le 2b$, we have
\begin{equation}
|\gamma'(s)|_{h(k)}\ge C^{-1}|\gamma'(b-s)|_{h(k)}
\end{equation}
where we have used the fact that
$|\gamma'(s)|_{h(k_0)}=|\gamma'(b-s)|_{h(k_0)}=1$. Hence
\begin{equation}
\ell_k(\gamma|_{(b\le s\le 2b)})\ge C^{-1}\ell_k(\gamma|_{(0\le s\le  b)})
\end{equation}
and
\begin{equation}
\begin{split}
\ell_k(\gamma|_{(0\le s\le 2b)})&\ge \lf(1+C^{-1}\ri)\ell_k(\gamma|_{(0\le s\le
b)})\\
&\ge \lf(1+C^{-1}\ri)b\\
&\ge \lf(1+C^{-1}\ri)(L-\epsilon).
\end{split}
\end{equation}
Let $k\to \infty$, we have
\begin{equation}
\ell_H(\gamma|_{(0\le s\le 2b)})\ge \lf(1+C^{-1}\ri)(L-\epsilon).
\end{equation}
Since $C$ does not depend on $\epsilon$, if we let $\epsilon\to0$, we have
\begin{equation}
\ell_H(\gamma|_{(0\le s\le 2b)})\ge \lf(1+C^{-1}\ri) L>L.
\end{equation}
This contradicts the definition of $L$.
\end{proof}
\begin{proof} (Sufficient part of Theorem \ref{mainthm} ): The first part of the conclusion follows from  Lemma \ref{existenceoflimit}  and Lemma
\ref{completeness}. In particular,   $H$ is a complete flat \K metric on $M$ and
thus $M$ is biholomorphic to a quotient of $\ce ^n$ by a group of biholomorphic
isometries.  But by Lemma  \ref{fixedpoint}  we know that $M$ is diffeomorphic to
$\re ^{2n}$.  Thus we must have $M$ biholomorphic to $\ce ^n$.
\end{proof}

\section{Convergence of K\"ahler-Ricci flows}

In this section we study a general solution to the K\"ahler-Ricci flow focusing
on Shi's program \cite{Sh} for the uniformization conjecture of Greene-Wu-Siu-Yau.  We will study the
\KRF equation

\begin{equation}\label{krf1}
\frac{\p }{\p t}g_\abb(x,t)=-R_\abb(x,t).
\end{equation}
More precisely, we are interested in the following situation. Let $(M^n,g_\abb)$
be a complete noncompact K\"ahler manifold with bounded nonnegative holomorphic
bisectional curvature such that the scalar curvature $R$ satisfies:
\begin{equation}\label{quaddecay}
\frac{1}{V_x(r)}\int_{B_x(r)}R \le \frac{C}{1+r^2}
\end{equation}
for some constant $C$ for all $x\in M$ and for all $r$. By \cite{Sh, Sh2, NT}, we
have the following:

\begin{thm}\label{shitamthm} Let $(M^n,g)$ be as above. Then the \KRF (\ref{krf1}) has long time solution with initial value $g_\abb(x,0)=g_\abb(x)$. Moreover, the following are true:
\begin{enumerate}
\item for any $t\ge 0$, $g(x,t)$ is K\"ahler with nonnegative holomorphic
bisectional curvature; \item for any $T>0$, there exists a constant $C_1>0$ such
that
$$
C_1^{-1}g(x,0)\le g(x,t)\le C_1g(x,0)
$$
for all $x\in M$ and for all $0\le t\le T$; \item for any integer $m\ge 0$, there
is a constant $C_2$ depending only on $m$ and the initial metric such that
$$
||\nabla^m Rm||^2(x,t)\le \frac{C_2}{(1+t)^{2+m}},
$$
for all $x\in M$ and for all $t$ if $m=0$ and for all $t\ge1$ if $m\ge1$, where $\nabla $ is the covariant derivative
with respect to $g(t)$ and the norm is also taken in $g(t)$.
\end{enumerate}
\end{thm}

 For the rest of the paper, we will always assume the conditions of Theorem \ref{shitamthm}.
 For any $T\ge 0$, define
$$
Q(x,t;T)=\lf(1+g^{\a\bd}(x,t)g^{\gamma\bb}(x,t)  g^{\xbz}(x,T)\tn_\x g_\abb(x,t)
\tn_\bz g_\gbd(x,t)\ri)^\frac12,
$$  where $t\ge T$,   and $\tn$ is the derivatives with respect to  $ g(T)$.  In
\cite{Sh}, a bound on $Q$ was derived in order to prove the existence of a
rescaled limit metric $h$ on $M$.  However, the derivation of this bound seems to
be incorrect. In particular the formula of $\frac{\p Q}{\p t}$ on p. 156 in \cite{Sh} is not correct.   Moreover, the proof of the completeness of $h$ is absent in \cite{Sh}.
In the first part of this section, we will  prove that the limit metric is complete under the assumption that a bound on $Q$ exists.

Let $p\in M$ be a fixed point and let $B_t(R)$ denote the geodesic ball  of
radius $R$ in $g(t)$ with center at $p$. Let $\mathbf v_p\in T^{1,0}(M)$ be a
fixed vector with length 1 in $g(0)$. As before, the norm of a vector in $g(t)$
is denoted by $|\mathbf v|_t$. We want to prove that:

\begin{thm}\label{firstway} Same assumptions as in Theorem \ref{shitamthm}. Moreover, suppose $M$ has positive holomorphic bisectional curvature and suppose for any $R>0$, there is a constant $C$ such that
\begin{equation}\label{qbound}
Q(x,t;T)\le C
\end{equation}
for all $T\ge0$, for all $x\in B_T(R)$ and for all $t\ge T$. Then there exists a
sequence $t_k\to\infty$ such that the metrics $\frac{1}{|\mathbf
v_p|^2_{t_k}}g(t_k)$ converge uniformly in $C^\infty$ topology to a complete
K\"ahler flat metric  on $M$. In particular, the universal covering space of $M$
is biholomorphic to $\mathbb C^n$.
\end{thm}
The crucial point is Lemma 3.3, which is also true under the assumptions of the
theorem.

\begin{lem}\label{uniformshrink1} With the same assumptions as in Theorem \ref{firstway}, let $R>0$ and   $T\ge0$. Then there exists a constant $C_R>0$ which is independent of $T$ with the following property:  For any $q\in B_{T}(R)$, $\mathbf w_p\in T^{1,0}(M)$, $\mathbf w_q\in T^{1,0}(M)$ with $|\mathbf w_p|_{T}=|\mathbf w_q|_{T}$,
$$
C_R^{-1}\le \frac{|\mathbf w_p|_{t}}{|\mathbf w_q|_t}\le C_R
$$
for all $t\ge T$. The constant $C_R$ is also independent of $q$, $\mathbf w_p$ $q$
and $\mathbf w_q$.
\end{lem}

\begin{proof}
This was basically proved in \cite{Sh}. Let $q$, $\mathbf w_p$, and $\mathbf w_q$
as in the lemma. Let $\gamma$ be a minimal geodesic with respect to $g(T)$ from
$q$ to $p$ parametrized by arc length and with length $\ell\le R$. Parallel
translate $\mathbf w_q$ along $\gamma$ with respect to $g(T)$ to obtain a vector
field $\mathbf w(s)$ on $\gamma$ such that $\mathbf w(0)=\mathbf w_q$. At any
point $s\in [0,\ell]$. Let $\tn$ be the covariant derivatives with respect to
$g(T)$. For any $t\ge T$, and for any $s$, choose an unitary frame near
$\gamma(s)$ such that $g_{\abb}(\gamma(s),T)=\delta_{\a\beta}$ and
$g_{\abb}(\gamma(s),t)=\lambda_\a\delta_{\a\beta}$. In the following, we write
$g=g(t)$ and $\wt g=g(T)$. Then \begin{equation}\label{gradg}
\begin{split}
\lf|\tn_{\mathbf \gamma'}\lf(g_\abb w^\a w^{\bb}\ri)\ri|^2&= \lf(\tn_{\mathbf
\gamma'}g_\abb\ri)
\lf(\tn_{ \mathbf \gamma' }g_\gbd\ri)w^\a w^{\bb}w^\gamma w^{\bar\delta}\\
&\le \lf(g_\abb w^\a w^{\bb}\ri)^2\sum_{\a,\beta}
\lambda_\a^{-1}\lambda_\beta^{-1}|\tn_{\mathbf u} g_\abb|^2
\end{split}
\end{equation} where we have used the facts that  $g_\abb w^\a
w^{\bb}=\sum_{\a}\lambda_\a|w^\a|^2$ and that $\mathbf w$ is parallel with
respect to $g(T)$, and the Schwarz inequality. On the other hand,
\begin{equation}\label{q2}
\begin{split}
Q^2&(\gamma(s),t;T)\\
&\ge g^{\a\bd} (\gamma(s),t)g^{\gamma\bb}(\gamma(s),t)g^{\xbz}(\gamma(s),T)\tn_\xi g_\abb(\gamma(s),t)\tn_{\bz} g_\gbd(\gamma(s),t)\\
&=\sum_{\a,\beta,\xi}\lambda_\a^{-1}\lambda_\beta^{-1}|\tn_{\xi} g_\abb|^2.
\end{split}
\end{equation}
 Combining (\ref{qbound}), (\ref{gradg}), (\ref{q2}) and the fact that $|\gamma'|_T=1$, we have
$$
\bigg|\tn_{\gamma'} \lf(g_\abb(\gamma(s),t_k) w^\a w^{\bb}\ri)\bigg|^2
 \le C_1 \lf(g_\abb w^\a w^{\bb}\ri)^2
 $$
for some constant $C_1$ which is independent of $t$, $T$, $q$ $\mathbf w_p$, and
$\mathbf w_q$. Hence and
$$
\bigg|\frac{\p}{\p s} \lf(g_\abb(\gamma(s),t_k) w^\a w^{\bb}\ri)\bigg|^2
 \le C_1 \lf(g_\abb v^\a v^{\bb}\ri)^2.
 $$
Integrating from $s=0$ to $s=\ell$, we have \begin{equation}\label{1}
\lf|\log\frac{g_\abb(t) v^\a(\ell) w^{\bb}(\ell)}{g_\abb(t) w^\a(0) w^{\bb}(0)}
\ri|\le C_3 \ell\le C_3R \end{equation} for some constant $C_3$ independent of
$t$, $T$, $q$ and $\mathbf w_q$ and $\mathbf w_p$. In particular, if we take
$t=0$, using the fact that the holomorphic bisectional curvature of $g(x,0)$ is
positive and hence the holonomy group is transitive \cite{simons}, we may prove
as in \cite{Sh} that \begin{equation}\label{2} \lf|\log\frac{g_\abb(t) u_1^\a
u_1^{\bb} }{g_\abb(t) u_2^\a u_2^{\bb} } \ri|\le C_4 \end{equation} for some
constant $C_4$, for all $\mathbf u_1, \ \mathbf u_2\in T_p^{1,0}(M)$ such that
$g_\abb(0)u_1^\a u_1^{\bb}=g_\abb(0)u_2^\a u_2^{\bb} $. Now if $\mathbf u_1$,
$\mathbf u_2$ are such that $g_\abb(T)u_1^\a  u_1^{\bb} =g_\abb(T)u_2^\a
u_2^{\bb}$, then by (\ref{2}), we have
$$
\frac{g_\abb(t) u_1^\a  u_1^{\bb} }{g_\abb(t) u_2^\a  u_2^{\bb} } \frac{g_\abb(0)
u_2^\a  u_2^{\bb} }  {g_\abb(0) u_1^\a  u_1^{\bb} }\le \exp(C_4)
$$
and
$$
\frac{g_\abb(T) u_2^\a  u_2^{\bb} }{g_\abb(T) u_1^\a  u_1^{\bb} } \frac
{g_\abb(0) u_1^\a  u_1^{\bb} }{g_\abb(0) u_2^\a  w_2^{\bb} }\le\exp( C_4)
$$
and hence we have
$$\frac{g_\abb(t) u_1^\a  u_1^{\bb} }{g_\abb(t) u_2^\a  u_2^{\bb} } \le \exp(2C_4)$$
for all $t\ge T$. Combining this to (\ref{1}), using the fact that
\begin{equation}
\begin{split}
g_\abb(T)w_p^\a w_p^\bb&=g_\abb(T)w_q^\a w_q^\bb\\
&=g_\abb(T)w^\a(\ell) w^\bb(\ell) \end{split} \end{equation} the lemma is proved.
\end{proof}

\begin{proof}(Theorem \ref{firstway}) Let $h(t)=\frac{1}{|\mathbf
v_p|^2}g(t)$.  By the proof of completeness in Lemma \ref{completeness}, because of
Lemma \ref{uniformshrink1} and Theorem \ref{shitamthm} it is sufficient to prove the
existence a limit for $h(t_k)$.  For this it is sufficient to show that in a fixed
coordinate neighborhood $U\subset B_0(R)$ the Christoffel symbols of $h(t)$ are
uniformly bounded. This can be proved as in Lemma \ref{existenceoflimit}. In this
case, using Theorem \ref{shitamthm}(3), Lemma \ref{uniformshrink1}, and the fact
that $g(t)$ is nonincreasing, we can conclude as in the proof of Lemma
\ref{existenceoflimit} that
\begin{equation}
\frac{\p}{\p t}||A||_0^2 \le C_1(1+t)^{-\frac 32}||A||_0
\end{equation}
where $A$ is defined as in Lemma \ref{existenceoflimit}, which is the difference
between the Christoffel symbols of $g(t)$ and $g(0)$ and $||A||_0$ is the norm of $A$ in $g(0)$. Here $C_1$ is a constant
depending only on $g(0)$, $R$ and the constant $C$ in the assumption
(\ref{qbound}) in the theorem. From this it is easy to see that $||A||_0$ is
uniformly bounded in $U\times[0,\infty)$. Hence the theorem is true.
\end{proof}

In the second part of this section, we will prove the following:

\begin{thm}\label{secondway} There exists a constant $C(n)$ depending only on $n$ such that if $M^n$ is a complete noncompact K\"ahler manifold satisfying the conditions in \ref{shitamthm} and the following:
\begin{enumerate}
\item[(i)]
 $$
\frac1{V_x(r)}\int_{B_x(r)}R\le \frac{C(n)}{1+r^2}
$$
for all $x\in M$ and for all $r>0$; and
\item[(ii)]   there exist  a point $p \in M$ and a sequence $t_k\to\infty$ such that   $\frac{1}{|\mathbf v_p|^2_{t_k}}g(p,t_k)$ are uniformly equivalent to $g(p,0)$, where $\mathbf v_p$ is a fixed vector in $T_p^{1,0}(M)$ with $|\mathbf v_p|_0=1$.
\end{enumerate}  Then the metrics $\frac{1}{|\mathbf v_p|^2_{t_k}}g(x,t_k)$
subconverge uniformly in the $C^\infty$ topology in compact sets to a complete K\"ahler flat metric
on $M$. In particular, the universal covering space of $M$ is biholomorphic to
$\mathbb C^n$.
\end{thm}

In order to prove the theorem, we need several lemma.

\begin{lem}\label{NiTam} Let $(M^n,  g)$ be a complete noncompact K\"ahler manifold with
nonnegative and bounded holomorphic bisectional curvature. Suppose there exists a
constant $a>0$ such that
\begin{equation}
\frac{1}{V_x(r)}\int_{B_x(r)}R\le \frac{a}{1+r^2}
\end{equation}
for all $x\in M$ for all $r$. Let $g_\abb(x,t)$ be the long time solution of (\ref{krf1}).
Then there exist constants $C_1$ depending only on $n$ and $C_2$ depending only
on $a$ and $n$ such that
\begin{equation}\label{NT0}
\int_0^t R(x,\tau)d\tau\le aC_1\log(1+t)+C_2
\end{equation}
for all $x\in M$ and for all $t$, where $R(x,t)$ is the scalar curvature of
$g(t)$ at $x$.
\end{lem}
\begin{proof}   For fixed $t$, the scalar curvature  $R(x,\tau)$ of $g(\tau)$ is uniformly bounded on $M\times[0,t]$. Let
\begin{equation}
\mathfrak M(t)=\max_{x\in M}\int_0^t R(x,\tau)d\tau.
\end{equation}
By \cite[Corollary2.1]{NT},  there exist positive constants $C_3$ and $C_4$
depending only on $n$ such that if $r^2=C_4t(1+\mathfrak M(t))$, then
\begin{equation}\label{NT1}
\begin{split}
 \mathfrak M(t)&\le C_3\int_0^{r}\frac{as}{1+s^2}ds\\
&= \frac{aC_3}{2}\log\lf(1+C_4t\lf(1+\mathfrak M(t)\ri)\ri)\\
&\le \frac{aC_3}2\lf[\log(1+C_4)+\log(1+t)+\log (1+\mathfrak M(t))\ri].
\end{split}
\end{equation}
Suppose $\mathfrak M(t)\ge aC_3$, then
\begin{equation}
\log(1+\mathfrak M(t))\le \frac1{aC_3}\mathfrak M(t)+\log(1+aC_3).
\end{equation}
By (\ref{NT1}), we have
\begin{equation}
\mathfrak M(t)\le aC_3\lf[\log(1+C_4)+\log(1+t)+\log(1+aC_3)\ri].
\end{equation}
Hence we have
\begin{equation}
\mathfrak M(t)\le aC_3\lf[\log(1+C_4) +\log(1+t)+\log(1+aC_3)+1\ri].
\end{equation}
From this  the lemma follows.
\end{proof}

\begin{lem}\label{Lemma0} Let $M^n$ be a complete noncompact K\"ahler manifold satisfying the conditions in Theorem \ref{shitamthm} and let $g(t)$ be the solution in (\ref{krf1}). Then for any $R>0$ there is a constant $C_R$ such that for any $T\ge
1$, $q\in B_T(R)$ and any $t\ge T$,
\begin{equation}
Rc_{\max}(q,t)-Rc_{\min}(q,t)\le Rc_{\max}(p,t)-Rc_{\min}(p,t)+C_R(1+t)^{-\frac
32}.
\end{equation}
Here as before, $p$ is a fixed point and $B_T(R)$ is the geodesic ball of radius
$R$ with center at $p$ in $g(T)$.
\end{lem}
\begin{proof} For $t\ge T$, et $\mathbf v_q$ and $\mathbf w_q$ in $T^{1,0}(M)$ such
that $|\mathbf v_q|_t=|\mathbf w_q|_t=1$ and $R_\abb(q,t)v_q^\a
v_q^\bb=Rc_{\max}(q,t)$, $R_\abb(q,t)w_q^\a w_q^\bb=Rc_{\min}(q,t)$. Let
$\gamma(s)$ be a minimal geodesic from $p$ to $q$ in $g(t)$ with length $\ell$
which is no greater than $R$ because $B_T(R)\subset B_t(R)$. Let $\mathbf v(s)$
and $\mathbf w(s)$ be parallel vector fields along $\gamma$ in $g(t)$ so that
$\mathbf v(\ell)=\mathbf v_q$ and $\mathbf w(\ell)=\mathbf w_q$. Then
\begin{equation}
\begin{split}
 &Rc_{\max}(q,t)-Rc_{\min}(q,t)\\
&=R_\abb(q,t)v_q^\a v_q^\bb- R_\abb(q,t)w_q^\a w_q^\bb\\
&=R_\abb(q,t)v(\ell)^\a v(\ell)^\bb- R_\abb(q,t)w(\ell)^\a w(\ell)^\bb\\
&=R_\abb(p,t)v(0)^\a v(0)^\bb- R_\abb(p,t)w(0)^\a w(0)^\bb\\
&\quad+\int_0^\ell\frac{\p}{\p s}\lf[ R_\abb(\g(s),t)v(s)^\a v(s)^\bb- R_\abb(\g(s),t)w(s)^\a w(s)^\bb \ri]ds\\
&=R_\abb(p,t)v(0)^\a v(0)^\bb- R_\abb(p,t)w(0)^\a w(0)^\bb\\
&\quad+\int_0^\ell \lf(\nabla_{\g'(s)}R_\abb(\g(s),t)\ri)v(s)^\a v(s)^\bb ds\\
&\quad-\int_0^\ell \lf( \nabla_{\g'(s)}R_\abb(\g(s),t)\ri) w(s)^\a w(s)^\bb ds\\
&\le Rc_{\max}(p,t)-Rc_{\min}(p,t)+CR(1+t)^{-\frac 32}
\end{split}
\end{equation}
where $C$ is a constant depending only on $g(0)$, where we have used Theorem \ref{shitamthm}. This completes the proof of the
lemma.
\end{proof}

\begin{lem}\label{lem1} With the same assumptions as in Lemma \ref{NiTam}, suppose $\epsilon=aC_1<1$ in (\ref{NT0}). For any $R>0$, $t\ge T\ge1$, $q\in B_T(R)$,
$\mathbf v_q,\ \mathbf w_q\in T^{1,0}(M)$ such that $|\mathbf v_q|_T=|\mathbf
w_q|_T$, we have
\begin{equation}
\frac{|\mathbf v_q|_t}{|\mathbf w_q|_t}\le C_R\lf(1+ t\ri)^{\frac12\epsilon}.
\end{equation}
where $C_R$ is a constant independent of $T,\ t, \ q,\ \mathbf v_q,\ \mathbf
w_q$.
\end{lem}
\begin{proof} By the proof of Lemma \ref{Lemma0} we have for $q\in B_T(R)$ and
$t\ge T$,
\begin{equation}
Rc_{\max}(q,t)-Rc_{\min}(q,t)\le Rc_{\max}(p,t)-Rc_{\min}(p,t)+C_1(1+t)^{-\frac
32},
\end{equation}
where $C_1$ is a constant independent of $T,\ t, \ q$.
Hence for $t\ge T$,
\begin{equation}
\begin{split}
\frac{\p}{\p t}\log \frac{|\mathbf v_q|_t}{|\mathbf w_q|_t}&\le \frac12\lf(-\frac{R_\abb(q,t)v_q^\a v_q^\bb}{g_\abb(q,t)v_q^\a v_q^\bb}+\frac{R_\abb(q,t)w_q^\a w_q^\bb}{g_\abb(q,t)w_q^\a w_q^\bb}\ri)\\
&\le \frac12\lf(Rc_{\max}(p,t)-Rc_{\min}(p,t)+C_R(1+t)^{-\frac 32}\ri)
\end{split}
\end{equation}
Integrating from $T$ to $t$ and using (\ref{NT0}), we have
\begin{equation}
\log \frac{|\mathbf v_q|_t}{|\mathbf w_q|_t}\le
\frac12\epsilon\log(1+t)+C_2,
\end{equation}
where $C_2$ is a constant independent of $T,\ t, \ q,\ \mathbf v_q,\ \mathbf
w_q$.
 Hence the lemma is true.
\end{proof}

As before, let $A^\gamma_{\a\b}=\Gamma^\gamma_{\a\b}-\wt\Gamma^\gamma_{\a\b}$,
where $\Gamma^\gamma_{\a\b}$ and $\wt\Gamma^\gamma_{\a\b}$ are Christoffel
symbols of $g(t)$ and $\wt g=g(T)$ respectively. Consider the norm of $A$ in
$g(T)$. Namely
\begin{equation}
||A||_T^2=\wt g_{\gamma\bar\tau}\wt g^{\a\bd} \wt g^{\xbz}A_{\alpha\x}^\gamma
A_{\bd\bz}^{\bar\tau}.
\end{equation}
\begin{lem}\label{lem3} With the same assumptions as in Lemma \ref{NiTam}, suppose $\epsilon=aC_1<1$ in (\ref{NT0}).  Then for any $R>0$, there is a constant $C_R$
such that for any $T\ge1$ and for any $t\ge T$,
\begin{equation}
||A||_T\le C_R
\end{equation}
in $B_T(R)$.
\end{lem}
\begin{proof}  As in (\ref{christoffel}), we have
\begin{equation}\label{estA}
\frac{\p}{\p t}||A||^2_T=-\wt g_{\gamma\bd}\wt g^{\abb} \wt
g^{\x\bz}\lf[g^{\gamma\bar\sigma}\nabla_\a R_{\x\bar\sigma}\ol{A_{\b\z}^\d}
+A_{\a\x}^\gamma \ol{g^{\d\bar\sigma}\nabla_\beta R_{\z\bar\sigma}}\ri]
\end{equation}
Choose  coordinates  such that $\wt g_\abb=\delta_{\a\b}$,
$g_\abb=\lambda_\a\delta_{\a\b}$. Then in $B_T(R)$,
\begin{equation}
\begin{split} \lf|\wt g_{\gamma\bd}\wt g^{\abb} \wt g^{\x\bz}  g^{\gamma\bar\sigma}\nabla_\a
R_{\x\bar\sigma}\ol{A_{\b\z}^\d}\ri|
&\le C(n) \sum_{\gamma,\a,\xi}\lambda_\g^{-1}\lf|\nabla_\a R_{\x\bar\gamma}\ri|\,||A||_T\\
&=C(n)\sum_{\gamma,\a,\xi}\lambda_{\a}^{\frac12}  \lambda_{\x}^\frac12 \lambda_\gamma^{-\frac12}
\lambda_\g^{-\frac12}\lambda_{\x}^{-\frac12}\lambda_{\a}^{-\frac12}\lf|\nabla_\a R_{\x\bar\gamma}\ri|\,||A||_T\\
&\le C_2 \lf(1+t\ri)^{\frac12\epsilon}\sum_{\gamma,\a,\xi} \lambda_\gamma^{-\frac12}\lambda_{\x}^{-\frac12}
\lambda_{\a}^{-\frac12} \lf|\nabla_\a R_{\x\bar\gamma}\ri|\,||A||_T\\
&\le C_3(1+t)^{-\frac 32+\frac12\epsilon}||A||_T
\end{split}
\end{equation}
for some constants $C_2$, $C_3$ independent of $t$ and $T$, where we have used
Lemma \ref{lem1}, the fact that $\lambda_\a\le1$ and the estimates for $||\nabla Rc||$.
Combining this with (\ref{estA}),  since $\epsilon<1$, and $||A||_T(T)=0$, it is easy to
see that the lemma is true.
\end{proof}

\begin{lem}\label{lem4} With the same assumptions as in Lemma \ref{NiTam}, suppose $\epsilon=aC_1<1$ in (\ref{NT0}). For any $R>0$ there is a constant $C_R$
such that if $t\ge T\ge1$, then
\begin{equation}\label{4.1}
C_R^{-1}\le \frac{g^\abb(x, T)g_\abb(x,t)}{g^\abb(p, T)g_\abb(p,t)}\le C_R
\end{equation}
and
\begin{equation}\label{4.2}
C_R^{-1}\le \frac{g^\abb(x, t)g_\abb(x,T)}{g^\abb(p, t)g_\abb(p,T)}\le C_R
\end{equation}
 for $x\in B_T(R)$.
\end{lem}
\begin{proof} We only prove (\ref{4.1}) as the proof of (\ref{4.2}) is similar.
We want to estimate $\lf|\tn \log \lf[g^\abb(x,
T)g_\abb(x,t)\ri]\ri|$ in $B_T(R)$, where $\tn$ is the covariant derivative of
$g(T)$. At a point, choose a normal coordinates so that $ g_\abb(T)=\delta_\abb$
and $g_\abb(t)=\lambda_\a\delta_{\a\b}$. Then
\begin{equation}
\begin{split}
\frac{\p}{\p \xi}\log \lf[g^\abb(x, T)g_\abb(x,t)\ri]&=\frac{g^\abb(x, T)\frac{\p}{\p \xi}g_\abb(x,t)}{g^\abb(x, T)g_\abb(x,t)}\\
&=\frac{g^\abb(x, T)\tn_\xi g_\abb(x,t)}{g^\abb(x, T)g_\abb(x,t)}\\
&=\frac{g^\abb(x, T)A^\tau_{\a\xi}g_{\tau\bb}(x,t)}{g^\abb(x, T)g_\abb(x,t)}.
\end{split}
\end{equation}
 Hence for $x$ in $B_T(R)$,
\begin{equation}
\begin{split} \lf|\frac{\p}{\p \xi}\log \lf[g^\abb(x, T)g_\abb(x,t)\ri]\ri|&\le
\frac{\sum_{\a}|A^\a_{\a\xi}|\lambda_\a}{\sum_\a\lambda_\a}\\&\le C_R,
\end{split}
\end{equation}
 where we have used Lemma \ref{lem3}. Hence
\begin{equation}\label{eq4.1}
\lf|\tn \log \lf[g^\abb(x, T)g_\abb(x,t)\ri]\ri|\le C_R
\end{equation}
in $B_T(R)$. Integrating (\ref{eq4.1}) along a minimal geodesic in $g(T)$ from
$x$ to $p$, (\ref{4.1}) follows.
\end{proof}

\begin{lem}\label{lem5} With the same assumptions as in Lemma \ref{NiTam}, suppose $\epsilon=aC_1<1$ in (\ref{NT0})  and suppose there exist $t_k\to\infty$, $t_k\ge 1$, such that
$\frac{1}{|\mathbf v_p|^2_{t_k}}g(p,t_k)$ are uniformly equivalent to $g(p,0)$, where $\mathbf v_p$ is a fixed vector in $T_p^{1,0}(M)$ with $|\mathbf v_p|_0=1$. Then for any $R>0$ there is a constant $C_R$ independent of $k$ and $k_0$ such that
$$
\frac{|\mathbf u_q|_{t_k}}{|\mathbf w_q|_{t_k}}\le C_R
$$
for all $q\in B_{t_{k_0}}(R)$, $k\ge k_0$ and $\mathbf u_q$, $\mathbf w_q\in
T_q^{1,0}(M)$ with $|\mathbf u_q|_{t_{k_0}}=|\mathbf w_q|_{t_{k_0}}$.
\end{lem}
\begin{proof} By the assumption, there is a constant $C>0$ independent of $k$ and $k_0$ such that
$$
g^{\abb}(p,t_k)g_{\abb}(p,t_{k_0})g^{\gamma\bd}(p,t_{k_0})g_{\gamma\bd}(p,t_k)\le C.
$$
From this and Lemma \ref{lem4}, the result follows.
\end{proof}
\begin{proof}(Theorem \ref{secondway}) By Theorem 4.1, Lemmas \ref{lem3}, \ref{lem5}, one can proceed as in the proof of Theorem \ref{firstway} to conclude that Theorem \ref{secondway} is true.
\end{proof}
\bibliographystyle{amsplain}

\end{document}